\documentclass[12pt]{amsart}
\usepackage{graphicx}
\usepackage{amssymb}

\newtheorem{thm}{Theorem}

\newtheorem{cor}{Corollary}

\DeclareMathOperator{\Des}{Des}

\DeclareMathOperator{\Sol}{Sol}

\title[Quasisymmetric functions]{A note on three types of quasisymmetric functions}

\author[T. K. Petersen]{T. Kyle Petersen}

\begin{document}
\begin{abstract}
In the context of generating functions for $P$-partitions, we revisit three flavors of quasisymmetric functions: Gessel's quasisymmetric functions, Chow's type B quasisymmetric functions, and Poirier's signed quasisymmetric functions. In each case we use the inner coproduct to give a combinatorial description (counting pairs of permutations) to the multiplication in: Solomon's type A descent algebra, Solomon's type B descent algebra, and the Mantaci-Reutenauer algebra, respectively. The presentation is brief and elementary, our main results coming as consequences of $P$-partition theorems already in the literature.
\end{abstract}

\maketitle

\section{Quasisymmetric functions and Solomon's descent algebra}

The ring of quasisymmetric functions is well-known (see \cite{Stanley2}, ch. 7.19). Recall that a quasisymmetric function is a formal series \[Q(x_1, x_2, \ldots ) \in \mathbb{Z}[[x_1, x_2,\ldots ]] \] of bounded degree such that the coefficient of $x_{i_1}^{\alpha_1} x_{i_2}^{\alpha_2} \cdots x_{i_k}^{\alpha_k}$ is the same for all $i_1 < i_2 < \cdots < i_k$ and all compositions $\alpha = (\alpha_1, \alpha_2, \ldots, \alpha_k)$. Recall that a composition of $n$, written $\alpha \models n$, is an ordered tuple of positive integers $\alpha = (\alpha_1, \alpha_2, \ldots, \alpha_k)$ such that $|\alpha| = \alpha_1 + \alpha_2 + \cdots + \alpha_k = n$. In this case we say that $\alpha$ has $k$ parts, or $\# \alpha = k$. We can put a partial order on the set of all compositions of $n$ by reverse refinement. The covering relations are of the form \[ (\alpha_1, \ldots, \alpha_i + \alpha_{i+1}, \ldots, \alpha_k ) \prec (\alpha_1, \ldots, \alpha_i, \alpha_{i+1}, \ldots, \alpha_k).\] Let $\mathcal{Q}sym_n$ denote the set of all quasisymmetric functions homogeneous of degree $n$. The ring of quasisymmetric functions can be defined as $\mathcal{Q}sym := \bigoplus_{n \geq 0} \mathcal{Q}sym_n$, but our focus will be on the quasisymmetric functions of degree $n$, rather than the ring as a whole.

The most obvious basis for $\mathcal{Q}sym_n$ is the set of \emph{monomial} quasisymmetric functions, defined for any composition $\alpha = (\alpha_1, \alpha_2, \ldots, \alpha_k) \models n$,
\[ M_{\alpha} := \sum_{i_1 < i_2 < \cdots < i_k} x_{i_1}^{\alpha_1} x_{i_2}^{\alpha_2} \cdots x_{i_k}^{\alpha_k}.\]
We can form another natural basis with the \emph{fundamental} quasisymmetric functions, also indexed by compositions,
\[ F_{\alpha} := \sum_{ \alpha \preccurlyeq \beta } M_{\beta},\] since, by inclusion-exclusion we can express the $M_{\alpha}$ in terms of the $F_{\alpha}$:
\[ M_{\alpha} = \sum_{ \alpha \preccurlyeq \beta } (-1)^{\#\beta - \#\alpha} F_{\beta}.\]
As an example, \[ F_{(2,1)} = M_{(2,1)} + M_{(1,1,1)} = \sum_{i < j} x_i^2 x_j + \sum_{i < j < k} x_i x_j x_k  = \sum_{i \leq j < k} x_i x_j x_k .\]

Compositions can be used to encode descent classes of permutations in the following way. Recall that a \emph{descent} of a permutation $\pi \in \mathfrak{S}_n$ is a position $i$ such that $\pi_i > \pi_{i+1}$, and that an \emph{increasing run} of a permutation $\pi$ is a maximal subword of consecutive letters $\pi_{i+1} \pi_{i+2} \cdots \pi_{i+r}$ such that $\pi_{i+1} < \pi_{i+2} < \cdots < \pi_{i+r}$. By maximality, we have that if $\pi_{i+1} \pi_{i+2} \cdots \pi_{i+r}$ is an increasing run, then $i$ is a descent of $\pi$ (if $i\neq 0$), and $i+r$ is a descent of $\pi$ (if $i+r \neq n$). For any permutation $\pi \in \mathfrak{S}_n$ define the \emph{descent composition}, $C(\pi)$, to be the ordered tuple listing (from left to right) the lengths of the increasing runs of $\pi$. If $C(\pi) = (\alpha_1, \alpha_2, \ldots, \alpha_k)$, we can recover the descent set of $\pi$:
\[ \Des(\pi) = \{ \alpha_1, \alpha_1 + \alpha_2, \ldots, \alpha_1 + \alpha_2 + \cdots + \alpha_{k-1} \}.\] Since $C(\pi)$ and $\Des(\pi)$ have the same information, we will use them interchangeably. For example, the permutation $\pi = (3,4,5,2,6,1)$ has $C(\pi) = (3,2,1)$ and $\Des(\pi) = \{ 3, 5\}$.

Recall (\cite{Stanley1}, ch. 4.5) that a $P$-partition is an order-preserving map from a poset $P$ to some (countable) totally ordered set. To be precise, let $P$ be any labeled partially ordered set (with partial order $<_P$) and let $S$ be any totally ordered countable set. Then $f: P \to S$ is a $P$-partition if it satisfies the following conditions:
\begin{enumerate}

\item $f(i) \leq f(j)$ if $i <_P j$

\item $f(i) < f(j)$ if $i <_P j$ and $i > j$ (as labels)

\end{enumerate}

We let $\mathcal{A}(P)$ (or $\mathcal{A}(P;S)$ if we want to emphasize the image set) denote the set of all $P$-partitions, and encode this set in the generating function
\[ \Gamma(P) := \sum_{ f \in \mathcal{A}(P)} x_{f(1)} x_{f(2)} \cdots x_{f(n)},\]
where $n$ is the number of elements in $P$ (we will only consider finite posets). If we take $S$ to be the set of positive integers, then it should be clear that $\Gamma(P)$ is always going to be a quasisymmetric function of degree $n$. As an easy example, let $P$ be the poset defined by $3 >_P 2 <_P 1$. In this case we have \[ \Gamma(P) = \sum_{ f(3) \geq f(2) < f(1)} x_{f(1)} x_{f(2)} x_{f(3)}.\]

We can consider permutations to be labeled posets with total order $\pi_1 <_{\pi} \pi_2 <_{\pi} \cdots <_{\pi} \pi_n$. With this convention, we have
\begin{align*}
 \mathcal{A}(\pi) = \{ f: [n] \to S  \mid & f(\pi_1) \leq f(\pi_2) \leq \cdots \leq f(\pi_n) \\
 & \mbox{ and } k \in \Des(\pi) \Rightarrow f(\pi_k) < f(\pi_{k+1})\},
\end{align*}
and
\[ \Gamma(\pi) = \sum_{ \substack{ i_1 \leq i_2 \leq \cdots \leq i_n \\ k \in \Des(\pi) \Rightarrow i_k < i_{k+1}}} x_{i_1} x_{i_2} \cdots x_{i_n}.\]
It is not hard to verify that in fact we have
\[ \Gamma(\pi) = F_{C(\pi)},\] so that generating functions for the $P$-partitions of permutations $\pi \in \mathfrak{S}_n$ form a basis for $\mathcal{Q}sym_n$.

We have the following theorem related to $P$-partitions of permutations, due to Gessel \cite{Gessel}.
\begin{thm}\label{thm:partition}
As sets, we have the bijection
\[
\mathcal{A}(\pi;ST) \leftrightarrow \coprod_{\sigma\tau = \pi} \mathcal{A}(\tau;S) \oplus \mathcal{A}(\sigma;T),
\]
where $ST$ is the cartesian product of the sets $S$ and $T$ with the lexicographic ordering and $\coprod$ denotes the disjoint union.
\end{thm}

Let $X = \{ x_1, x_2, \ldots \}$ and $Y = \{ y_1, y_2, \ldots \}$ be two two sets of commuting indeterminates. Then we define the bipartite generating function, \[ \Gamma(\pi)(XY) = \sum_{\substack{ (i_1,j_1) \leq (i_2, j_2) \leq \cdots \leq(i_n,j_n) \\ k \in \Des(\pi) \Rightarrow (i_k,j_k) < (i_{k+1}, j_{k+1}) }} {\kern -10pt} x_{i_1}\cdots x_{i_n} y_{j_1} \cdots y_{j_n}. \] We will apply Theorem \ref{thm:partition} with $S = T = \mathbb{P}$, the positive integers.

\begin{cor}\label{cor:internal}
We have
\[ F_{C(\pi)}(XY) = \sum_{ \sigma\tau = \pi} F_{C(\tau)}(X) F_{C(\sigma)}(Y) .\]
\end{cor}

Following \cite{Gessel}, we can define a coalgebra $\mathcal{Q}sym_n^{*}$ in the following way. If $\pi$ is any permutation with $C(\pi) = \gamma$, let $a_{\alpha,\beta}^{\gamma}$ denote the number of pairs of permutations $(\sigma, \tau) \in \mathfrak{S}_n \times \mathfrak{S}_n$ with $C(\sigma) = \alpha$, $C(\tau) = \beta$, and $\sigma\tau = \pi$. Then Corollary \ref{cor:internal} defines a coproduct $\mathcal{Q}sym_n^* \to \mathcal{Q}sym_n^* \otimes \mathcal{Q}sym_n^*$: \[ F_{\gamma} \mapsto \sum_{\alpha,\beta \models n} a_{\alpha,\beta}^{\gamma} F_{\beta} \otimes F_{\alpha}.\] The dual space to $\mathcal{Q}sym_n^{*}$ is then $\mathcal{Q}sym_n$ equipped with multiplication \[ F_{\beta} \ast F_{\alpha} = \sum_{\gamma} a_{\alpha,\beta}^{\gamma\models} F_{\gamma}.\]

Let $\mathbb{Z}\mathfrak{S}_n$ denote the group algebra of the symmetric group. We can define its dual coalgebra $\mathbb{Z}\mathfrak{S}_n^*$ with comultiplication \[ \pi \mapsto \sum_{\sigma \tau = \pi} \tau \otimes \sigma.\] Then we have a surjective homomorphism of coalgebras $\varphi^*: \mathbb{Z}\mathfrak{S}_n^* \to \mathcal{Q}sym_n^*$ given by \[ \varphi^*(\pi) =  F_{C(\pi)}.\] The dualization of this map is then an injective homomorphism $\varphi: \mathcal{Q}sym_n \to \mathbb{Z}\mathfrak{S}_n$ with \[ \varphi(F_{\alpha}) = \sum_{C(\pi) = \alpha} \pi.\]

The image of $\varphi$ is then a subalgebra of the group algebra, with basis \[ u_{\alpha} := \sum_{C(\pi) = \alpha} \pi.\] This subalgebra is well-known as \emph{Solomon's descent algebra}, denoted $\Sol(A_{n-1})$. Corollary \ref{cor:internal} has then given a combinatorial description to multiplication in $\Sol(A_{n-1})$:
\[ u_{\beta}u_{\alpha} = \sum_{\gamma \models n} a_{\alpha, \beta}^{\gamma} u_{\gamma}.\]

The above arguments are due to Gessel \cite{Gessel}. We give them here in full detail for comparison with later sections, when we will outline a similar relationship between type Chow's B quasisymmetric functions and $\Sol(B_n)$, and between Poirier's \emph{signed} quasisymmetric functions and the Mantaci-Reutenauer algebra.

\section{Type B quasisymmetric functions and Solomon's descent algebra}

The type B quasisymmetric functions can be viewed as the natural objects related to type B $P$-partitions (see \cite{Chow}). Define the type B posets (with $2n+1$ elements) to be posets labeled distinctly by $\{-n, \ldots, -1, 0, 1, \ldots, n\}$ with the property that if $i <_P j$, then $-j <_P -i$. For example, $-2 >_P 1 <_P 0 <_P -1 >_P 2$ is a type B poset.

Let $P$ be any type B poset, and let $S = \{ s_0, s_1, \ldots \}$ be any countable totally ordered set with a minimal element $s_0$. Then a type B $P$-partition is any map $f: P \to \pm S $ such that
\begin{enumerate}

\item $f(i) \leq f(j)$ if $i <_P j$

\item $f(i) < f(j)$ if $i <_P j$ and $i > j$ (as labels)

\item $f(-i) = -f(i)$

\end{enumerate}
where $\pm S $ is the totally ordered set \[ \cdots < -s_2 < -s_1 < s_0 < s_1 < s_2 < \cdots \] If $S$ is the nonnegative integers, then $\pm S$ is the set of all integers.

The third property of type B $P$-partitions means that $f(0) = 0$ and the set $\{ f(i) \mid i = 1,2,\ldots,n\}$ determines the map $f$. We let $\mathcal{A}_B (P) = \mathcal{A}_B (P;\pm S)$ denote the set of all type B $P$-partitions, and define the generating function for type B $P$-partitions as
\[ \Gamma_B(P) := \sum_{f \in \mathcal{A}_B (P)} x_{|f(1)|} x_{|f(2)|} \cdots x_{|f(n)|}.\]

Signed permutations $\pi \in \mathfrak{B}_n$ are type B posets with total order \[ -\pi_n < \cdots < -\pi_1 < 0 < \pi_1 < \cdots < \pi_n. \] We then have
\begin{align*}
 \mathcal{A}_B(\pi) = \{ f: \pm [n] \to \pm S  \mid & f(-i) = -f(i), \\
 & 0 \leq f(\pi_1) \leq f(\pi_2) \leq \cdots \leq f(\pi_n), \\
 & \mbox{and } k \in \Des_B(\pi) \Rightarrow f(\pi_k) < f(\pi_{k+1})\},
\end{align*}
and
\[ \Gamma_B(\pi) = \sum_{ \substack{ 0\leq i_1 \leq i_2 \leq \cdots \leq i_n \\ k \in \Des(\pi) \Rightarrow i_k < i_{k+1}}} x_{i_1} x_{i_2} \cdots x_{i_n}.\]
Here, the type B descent set, $\Des_B(\pi)$, keeps track of the ordinary descents as well as a descent in position 0 if $\pi_1 < 0$. Notice that if $\pi_1 < 0$, then $f(\pi_1) > 0$, and $\Gamma_B(\pi)$ has no $x_0$ terms, as in \[ \Gamma_B((-3,2,-1) ) = \sum_{ 0 < i \leq j < k } x_i x_j x_k.\]

The possible presence of a descent in position zero is the crucial difference between type A and type B descent sets. Define a \emph{psuedo-composition} of $n$ to be an ordered tuple $\alpha = ( \alpha_1, \ldots, \alpha_k )$ with $\alpha_1 \geq 0$, and $\alpha_i > 0$ for $i > 1$, such that $\alpha_1 + \cdots + \alpha_k = n$. We write $\alpha \Vdash n$ to mean $\alpha$ is a psuedo-composition of $n$. Define the descent psuedo-composition $C(\pi)$ of a signed permutation $\pi$ be the lengths of its increasing runs as before, but now we have $\alpha_1 = 0$ if $\pi_1 < 0$. As with ordinary compositions, the partial order on pseudo-compositions of $n$ is given by reverse refinement. We can move back and forth between descent psuedo-compositions and descent sets in exactly the same way as for type A. If $C(\pi) = (\alpha_1, \ldots, \alpha_k)$, then we have \[ \Des_B (\pi) = \{ \alpha_1, \alpha_1 + \alpha_2, \ldots, \alpha_1 + \alpha_2 + \cdots + \alpha_{k-1} \}. \]

We will use psuedo-compositions of $n$ to index the type B quasisymmetric functions. Define $\mathcal{BQ}sym_n$ as the vector space of functions spanned by the \emph{type B monomial quasisymmetric functions}:
\[ M_{B,\alpha} := \sum_{ 0 < i_2 < \cdots < i_k } x_0^{\alpha_1} x_{i_2}^{\alpha_2} \cdots x_{i_k}^{\alpha_k}, \] where $\alpha = (\alpha_1,\alpha_2, \ldots, \alpha_k)$ is any psuedo-composition of $n$, or equivalently by the \emph{type B fundamental quasisymmetric functions}:
\[ F_{B,\alpha} := \sum_{ \alpha \preccurlyeq \beta } M_{B,\beta}.\] The space of all type B quasisymmetric functions is defined as the direct sum $\mathcal{BQ}sym := \bigoplus_{n\geq 0} \mathcal{BQ}sym_n$. By design, we have \[ \Gamma_B(\pi) = F_{B,C(\pi)}.\]

From Chow \cite{Chow} we have the following theorem and corollary.
\begin{thm}\label{thm:partitionB}
As sets, we have the bijection
\[
\mathcal{A}_B(\pi;ST) \leftrightarrow \coprod_{\sigma\tau = \pi} \mathcal{A}_B(\tau;S) \oplus \mathcal{A}_B(\sigma;T).
\]
\end{thm}

We take $S = T = \mathbb{Z}$ and we have the following.
\begin{cor}\label{cor:internalB}
We have
\[ F_{B,C(\pi)}(XY) = \sum_{ \sigma\tau = \pi} F_{B,C(\tau)}(X) F_{B,C(\sigma)}(Y) .\]
\end{cor}

The coalgebra structure on $\mathcal{BQ}sym_n$ works just the same as the type A case so we will omit some details. Corollary \ref{cor:internalB} gives us the coproduct \[ F_{B,\gamma} \mapsto \sum_{\alpha,\beta \Vdash n} b_{\alpha,\beta}^{\gamma} F_{B,\beta} \otimes F_{B,\alpha}, \] where for any $\pi$ such that $C(\pi) = \gamma$, $b_{\alpha,\beta}^{\gamma}$ is the number of pairs of signed permutations $(\sigma, \tau)$ such that $C(\sigma) = \alpha$, $C(\tau) = \beta$, and $\sigma\tau = \pi$. The dual algebra is isomorphic to $\Sol(B_n)$, where if $u_{\alpha}$ is the sum of all signed permutations with descent psuedo-composition $\alpha$, the multiplication given by
\[ u_{\beta}u_{\alpha} = \sum_{\gamma \Vdash n} b_{\alpha, \beta}^{\gamma} u_{\gamma}.\]

\section{Signed quasisymmetric functions and the Mantaci-Reutenauer algebra}

One thing to have noticed about the generating function for type B $P$-partitions is that we are losing a certain amount of information when we take absolute values on the subscripts. We can think of signed quasisymmetric functions as arising naturally by dropping this restriction.

For a type B poset $P$, define the \emph{signed generating function} for type B $P$-partitions to be
\[ \overline{\Gamma}(P) := \sum_{ f \in \mathcal{A}_B(P)} x_{f(1)} x_{f(2)} \cdots x_{f(n)},\]
where we will write
\[ x_i = \begin{cases} u_i & \mbox{ if } i < 0, \\ v_i & \mbox{ if } i \geq 0.
\end{cases} \]
In the case where $P$ is a signed permutation, we have
\[ \overline{\Gamma}(\pi) = \sum_{ \substack{ 0 \leq i_1 \leq i_2 \leq \cdots \leq i_n \\  s\in \Des_B(\pi) \Rightarrow i_s < i_{s+1} \\ \pi_s < 0 \Rightarrow x_{i_s} = u_{i_s} \\ \pi_s > 0 \Rightarrow x_{i_s} = v_{i_s}} } x_{i_1} x_{i_2} \cdots x_{i_n}, \] so that now we are keeping track of the set of minus signs of our signed permutation along with the descents. For example, \[ \overline{\Gamma}( (-3, 2, -1) ) = \sum_{ 0 < i \leq j < k } u_i v_j u_k.\]

To keep track of both the set of signs and the set of descents, we introduce the \emph{signed compositions} as used in \cite{BonaffeHohlweg}. A signed composition $\alpha$ of $n$, denoted $\alpha \Vvdash n$, is a tuple of nonzero integers $(\alpha_1, \ldots, \alpha_k)$ such that $|\alpha_1| + \cdots + |\alpha_k| = n$. For any signed permutation $\pi$ we will associate a signed composition $sC(\pi)$ by simply recording the length of increasing runs with constant sign, and then recording that sign. For example, if $\pi = ( -3, 4, 5,-6, -2, -7, 1)$, then $sC(\pi) = (-1, 2, -2, -1,1)$. The signed composition keeps track of both the set of signs and the set of descents of the permutation, as we demonstrate with an example. If $sC(\pi) = ( -3, 2, 1, -2, 1)$, then we know that $\pi$ is a permutation in $\mathfrak{S}_9$ such that $\pi_4, \pi_5, \pi_6$, and $ \pi_9$ are positive, whereas the rest are all negative. The descents of $\pi$ are in positions 5 and 6. Note that for any ordinary composition of $n$ with $k$ parts, there are $2^k$ signed compositions, leading us to conclude that there are \[ \sum_{k=1}^n \binom{n-1}{k-1} 2^k = 2\cdot3^{n-1} \] signed compositions of $n$.

We will use signed compositions to index the signed quasisymmetric functions (see \cite{Poirier}). For any signed composition $\alpha$, define the \emph{monomial signed quasisymmetric function}
\[ \overline{M}_{\alpha} := \sum_{ \substack{ i_1 < i_2 < \cdots < i_k \\ \alpha_r < 0 \Rightarrow x_{i_r} = u_{i_r} \\ \alpha_r > 0 \Rightarrow x_{i_r} = v_{i_r} } } x_{i_1}^{|\alpha_1|} x_{i_2}^{|\alpha_2|} \cdots x_{i_k}^{|\alpha_k|},\] and the \emph{fundamental signed quasisymmetric function}
\[ \overline{F}_{\alpha} := \sum_{ \alpha \preccurlyeq \beta} \overline{M}_{\beta}.\] By construction, we have \[ \overline{\Gamma}(\pi) = \overline{F}_{sC(\pi)}.\] Notice that if we set $u = v$, then our signed quasisymmetric functions become type B quasisymmetric functions.

Let $\mathcal{SQ}sym_n$ denote the span of the $\overline{M}_{\alpha}$ (or $\overline{F}_{\alpha}$), taken over all $\alpha \Vvdash n$. The space of all signed quasisymmetric functions, $\mathcal{SQ}sym := \bigoplus_{n\geq 0} \mathcal{SQ}sym_n$, is a graded ring whose $n$-th graded component has rank $2\cdot3^{n-1}$.

Theorem \ref{thm:partitionB} is a statement about splitting apart bipartite $P$-partitions, independent of how we choose to encode the information. So while Corollary \ref{cor:internalB} is one such way of encoding the information of Theorem \ref{thm:partitionB}, the following is another.

\begin{cor}\label{cor:internalS}
We have
\[ \overline{F}_{sC(\pi)}(XY) = \sum_{ \sigma\tau = \pi} \overline{F}_{sC(\tau)}(X) \overline{F}_{sC(\sigma)}(Y). \]
\end{cor}

We define the coalgebra $\mathcal{SQ}sym_n^*$ as we did in the earlier cases. Let $\pi \in \mathfrak{B}_n$ be any signed permutation with $sC(\pi) = \gamma$, and let $c_{\alpha,\beta}^{\gamma}$ be the number of pairs of permutations $(\sigma, \tau) \in \mathfrak{B}_n \times \mathfrak{B}_n$ with $sC(\sigma) = \alpha$, $sC(\tau) = \beta$, and $\sigma\tau = \pi$. Corollary \ref{cor:internalS} gives a coproduct $\mathcal{SQ}sym_n^* \to \mathcal{SQ}sym_n^* \otimes \mathcal{SQ}sym_n^*$: \[ \overline{F}_{\gamma} \mapsto \sum_{\alpha,\beta \Vvdash n} c_{\alpha,\beta}^{\gamma} \overline{F}_{\beta} \otimes \overline{F}_{\alpha}.\] Multiplication in the dual algebra $\mathcal{SQ}sym_n$, the signed quasisymmetric functions of degree $n$, is given by \[ \overline{F}_{\beta} \ast \overline{F}_{\alpha} = \sum_{\gamma \Vvdash n} c_{\alpha,\beta}^{\gamma} \overline{F}_{\gamma}.\]

The group algebra of the hyperoctahedral group, $\mathbb{Z}\mathfrak{B}_n$, has a dual coalgebra $\mathbb{Z}\mathfrak{B}_n^*$ with comultiplication given by the map \[ \pi \mapsto \sum_{\sigma \tau = \pi} \tau \otimes \sigma.\] The following is a surjective homomorphism of coalgebras $\psi^*: \mathbb{Z}\mathfrak{B}_n^* \to \mathcal{SQ}sym_n^*$ given by \[ \psi^*(\pi) =  \overline{F}_{sC(\pi)}.\] The dualization of this map is an injective homomorphism $\psi: \mathcal{SQ}sym_n \to \mathbb{Z}\mathfrak{B}_n$ with \[ \psi(\overline{F}_{\alpha}) = \sum_{sC(\pi) = \alpha} \pi.\]

The image of $\psi$ is then a subalgebra of $\mathbb{Z}\mathfrak{B}_n$ of dimension $2\cdot 3^{n-1}$, with basis \[ v_{\alpha} := \sum_{sC(\pi) = \alpha} \pi.\] This subalgebra is called the \emph{Mantaci-Reutenauer algebra}, with multiplication given explicitly by
\[ v_{\beta}v_{\alpha} = \sum_{\gamma \Vvdash n} c_{\alpha, \beta}^{\gamma} v_{\gamma}.\]

In closing, we remark that this same method is seen in \cite{Petersen2}, where Stembridge's enriched $P$-partitions \cite{Stembridge} are generalized and put to use to study peak algebras. Variations on the theme can also be found in \cite{Petersen1}.

\end{document}